\catcode`\@=11

\magnification=1200
\baselineskip=14pt

\pretolerance=500    \tolerance=1000 \brokenpenalty=5000

\catcode`\;=\active
\def;{\relax\ifhmode\ifdim\lastskip>\z@
\unskip\fi\kern.2em\fi\string;}

\overfullrule=0mm

\catcode`\!=\active
\def!{\relax\ifhmode\ifdim\lastskip>\z@
\unskip\fi\kern.2em\fi\string!}

\catcode`\?=\active
\def?{\relax\ifhmode\ifdim\lastskip>\z@
\unskip\fi\kern.2em\fi\string?}

\frenchspacing

\newif\ifpagetitre            \pagetitretrue
\newtoks\hautpagetitre        \hautpagetitre={\hfill \date }
\newtoks\baspagetitre         \baspagetitre={1}

\newtoks\auteurcourant        \auteurcourant={   }
\newtoks\titrecourant
\titrecourant={  }

\newtoks\hautpagegauche       \newtoks\hautpagedroite
\hautpagegauche={\hfill\sevenrm\the\auteurcourant\hfill}
\hautpagedroite={\hfill\sevenrm\the\titrecourant\hfill}

\newtoks\baspagegauche       \baspagegauche={\hfill\rm\folio\hfill}

\newtoks\baspagedroite       \baspagedroite={\hfill\rm\folio\hfill}

\headline={
\ifpagetitre\the\hautpagetitre
\global\pagetitrefalse
\else\ifodd\pageno\the\hautpagedroite
\else\the\hautpagegauche\fi\fi}

\footline={\ifpagetitre\the\baspagetitre
\global\pagetitrefalse
\else\ifodd\pageno\the\baspagedroite
\else\the\baspagegauche\fi\fi}

\def\date{\ {\the\day}\
\ifcase\month\or Janvier\or F\'evrier\or Mars\or Avril
\or Mai \or Juin\or Juillet\or Ao\^ut\or Septembre
\or Octobre\or Novembre\or D\'ecembre\fi\
{\the\year}}

\def\up#1{\raise 1ex\hbox{\sevenrm#1}}

\def\cqfd{\unskip\kern 6pt\penalty 500
\raise -2pt\hbox{\vrule\vbox to 10pt{\hrule width 4pt
\vfill\hrule}\vrule}\par\medskip}

\def\section#1{\vskip 7mm plus 20mm minus 1.5mm\penalty-50
\vskip 0mm plus -20mm minus 1.5mm\penalty-50
{\bf\noindent#1}\nobreak\smallskip}

\def\subsection#1{\medskip{\bf#1}\nobreak\smallskip}

\def\displaylinesno #1{\dspl@y\halign{
\hbox to\displaywidth{$\@lign\hfil\displaystyle##\hfil$}&
\llap{$##$}\crcr#1\crcr}}

\def\ldisplaylinesno #1{\dspl@y\halign{
\hbox to\displaywidth{$\@lign\hfil\displaystyle##\hfil$}&
\kern-\displaywidth\rlap{$##$}
\tabskip\displaywidth\crcr#1\crcr}}

\def\hfl#1#2{\smash{\mathop{\hbox to 12 mm{\rightarrowfill}}
\limits^{\scriptstyle#1}_{\scriptstyle#2}}}

\catcode`\@=12

\def\bN{{\bf N}}
\def\bP{{\bf P}}

\def\bR{{\bf R}}
\def\bZ{{\bf Z}}

\def\cO{{\cal O}}

\def\cS{{\cal S}}

\def\io{{\it o}}

\def\ux{{{\bf x}}}

\def\uy{{\bf y}}

\def\uz{{\bf z}}

\def\bZ{{\bf Z}}

\def\proof{\bigskip\noindent{\it Proof.}\ }

\def\and{\quad\hbox{and}\quad}


\def\M{\mathop{\rm M\kern 1pt}\nolimits}
\def\h{\mathop{\rm h\kern 1pt}\nolimits}

\def\romain#1{\uppercase\expandafter{\romannumeral #1}}

\def\card{\mathop{\rm Card\kern 1.3 pt}\nolimits}
\def\deg{\mathop{\rm deg\kern 1pt}\nolimits}
\def\det{\mathop{\rm det\kern 1pt}\nolimits}

\def\h{\mathop{\rm h\kern 1pt}\nolimits} \long\def\forget#1\endforget{}

\def\hw{{\hat{w}}}
\def\eps{{\epsilon}}
\def\al{{\alpha}}
\def\og{\leavevmode\raise.3ex\hbox{$\scriptscriptstyle
\langle\!\langle\,$}}
\def\fg{\leavevmode\raise.3ex\hbox{$\scriptscriptstyle
\!\rangle\!\rangle\,\,$}}
\def\thi{{\thinspace}}


\def\AdDa{1}
\def\ADQZ{2}
\def\ArRo{3}
\def\BaSc{4}
\def\RCBak{5}
\def\Bern{6}
\def\BuGr{7}
\def\BuJNT{8}
\def\BuAA{9}
\def\BuLiv{10}
\def\Cassa{11}
\def\DaScA{12}
\def\DaScB{13}
\def\DaScC{14}
\def\Dav{15}
\def\Kok{16}
\def\LauA{17}
\def\LauB{18}
\def\Mah{19}
\def\Que{20}
\def\RoSz{21}
\def\RoyA{22}
\def\RoyB{23}
\def\RoyC{24}
\def\RoyD{25}
\def\SprLiv{26}
\def\Wir{27}


\centerline{}

\vskip 4mm

\centerline{
\bf Exponents of Diophantine Approximation and Sturmian Continued Fractions}

\vskip 8mm
\centerline{Yann B{\sevenrm UGEAUD} \footnote{}{\rm
} \&
Michel L{\sevenrm AURENT} }

\vskip 9mm

\vskip 2mm
\noindent {\bf Abstract --} Let $\xi$ be a real number and let $n$ be
a positive integer.
We define four exponents of
Diophantine approximation, which complement the exponents  $w_n(\xi)$
and $w_n^*(\xi)$ defined
by Mahler and Koksma. We calculate their six values when $n=2$ and
$\xi$ is a real number
whose continued fraction expansion coincides with some
Sturmian sequence of positive integers, up to the initial terms.
   In particular, we obtain the exact exponent of
   approximation to such a continued fraction $\xi$ by quadratic surds.

\vskip15mm
\section{1. Introduction}

Mahler [\Mah] and Koksma [\Kok] have introduced two classifications
of real numbers $\xi$ in terms of their properties of approximation
by algebraic numbers.
  Keeping their notations, for every integer $n \ge 1$,
let us denote by $w_n(\xi)$ the supremum of the real numbers $w$ for which
the inequality
$$
0 < |P(\xi)| \le H(P)^{-w}
$$
is  satisfied for infinitely many polynomials
$P(X)$ with integer coefficients and degree at most $n$
(the height $H(P)$ of a polynomial $P(X)$ is the maximum of the moduli of
its coefficients).
In a similar way,
define $w_n^*(\xi)$ as the supremum of the real numbers $w^*$ for which 
the inequality
$$
0 < |\xi - \alpha| \le H(\alpha)^{- w^* -1}  \eqno (1)
$$
is  satisfied for infinitely many algebraic numbers
$\alpha$ of degree at most $n$
(the height $H(\alpha)$ of an algebraic number $\alpha$
is the height of its minimal polynomial over $\bZ$). The adjunction of
$-1$ in the exponent of the right hand side of (1) has the following 
motivation.
Assume that $|P(\xi)|$ is small. Then, the polynomial $P(X)$ has some root
  $\alpha$ with  $|\xi - \alpha| \ll |P(\xi) /
P'(\xi)|$; since $|P'(\xi)|$ should be {\it in  principle } roughly equal to
$H(P)$, we expect the estimation $|\xi - \alpha| \ll |P(\xi)|  \cdot
H(P)^{-1}$.

The behaviour of the sequences $(w_n (\xi))_{n \ge 1}$ and
$(w_n^*(\xi))_{n \ge 1}$ determines the localisation of $\xi$ in
  Mahler's and Koksma's classifications, respectively;
however, the precise estimation of $w_n (\xi)$ and $w_n^*(\xi)$
is usually extremely difficult.
The Dirichlet box principle (or, equivalently,
Minkowski's theorem) readily implies that $w_n(\xi) \ge n$ for
any positive integer $n$ and any transcendental real number $\xi$.
It is a longstanding problem to decide whether the same
result remains true for the functions $w_n^*$.

\proclaim Conjecture (Wirsing). For any integer $n \ge 1$
and any transcendental real number $\xi$ we have $w_n^*(\xi) \ge n$.

The early paper of Wirsing [\Wir] and the study of his
conjecture, which has been up to now confirmed only for $n=1$ 
(this is a simple application
of Dirichlet's box principle) and $n=2$ (by Davenport \& Schmidt [\DaScA]),
have motivated many works. In particular, Davenport \&
Schmidt [\DaScB] have investigated this question in a dual way,
via  simultaneous rational approximation to the successive powers of
$\xi$. Among other results, they proved that if the
real number $\xi$ is either 
transcendental, or algebraic of degree  $\ge 3$, then there
exist some  positive constant $c$, depending only upon  $\xi$, and
arbitrarily large real numbers X such that the inequalities
$$
0 < |x_0| < X, \quad |x_0 \xi - x_1| < c X^{-(\sqrt{5} - 1)/2}, \quad
|x_0 \xi^2 - x_2| < c X^{-(\sqrt{5} - 1)/2} \eqno (2)
$$
have no integer solution $x_0$, $x_1$, $x_2$.
This assertion is by no means surprising
since the  inequalities (2)
define a  convex body whose  volume tends to $0$ when $X$
tends to infinity, and  one should be tempted to believe that
the same conclusion holds if we  replace $(\sqrt{5} - 1)/2$
by any
real number larger than  $1/2 \ldots$, but that is false.
Indeed, when $\xi$ is a real number whose continued fraction expansion
reproduces the
structure  of the Fibonacci word, Roy proved
recently that there does exist an other $c$ such that $(2)$
    has an integer solution $x_0,x_1,x_2$ for any sufficiently large $X$
(Theorem 1.1 from [\RoyB]).  The critical exponent in $(2)$
is definitively $(\sqrt{5} - 1)/2$
for a Fibonacci continued fraction $\xi$.

The above point (2), as well as the article
[\ArRo], suggests to us to introduce other exponents which implicitely appear
in problems of  Diophantine Approximation by algebraic numbers or
algebraic integers. We define  four new exponents, which
naturally complement the  functions $w_n$ and $w_n^*$.

We first denote  by  $\lambda_n(\xi)$ the exponent
of simultaneous rational approximation to the numbers $\xi, \dots , 
\xi^n$. Thus,
$\lambda_n(\xi)$ denotes
the supremum of the real numbers $\lambda$ such that
the inequality
$$
\qquad  \max_{1 \le m \le n} \,
|x_0 \xi^m - x_m| \le  |x_0|^{-\lambda}
$$
has infinitely  many solutions in integers $x_0, \dots,x_n$
with $x_0\not= 0$. It will also be convenient to introduce
the inverse
$$
w'_n(\xi) = {1\over \lambda_n(\xi)},
$$
    which is equal to  $n$ for almost all real number $\xi$ (with 
respect to the Lebesgue measure).
The  three exponents $w_n$, $w^*_n$ and $\lambda_n$
have the common feature to be defined by the occurrence
of infinitely  many solutions for some set of Diophantine inequalities.
We attach to them three  functions defined
by a condition of uniform existence of solutions and
we decide to indicate this uniformity by the symbol $\, \hat{ }$ .

\proclaim
Definition 1.1.
Let $n \ge 1$ be an integer and let $\xi$ be a real number.
We denote by $\hw_n(\xi)$
the supremum of the real numbers $w$ such that, for
any sufficiently large real number $X$, the inequalities
$$
0 < |x_n \xi^n + \ldots + x_1 \xi + x_0| \le X^{-w}, \qquad  \max_{0
\le m \le n} \, |x_m| \le X,
$$
have a solution in integers $x_0, \ldots, x_n$.
We denote by  $\hw^*_n(\xi)$  the
  supremum of the real numbers $w$ such that, for
any sufficiently large real number $X$,
there exists an algebraic real number
  $\alpha$ with degree at most $n$ and satisfying
$$
0 < |\xi - \alpha| \le H(\alpha)^{-1} \, X^{-w} \and H(\alpha) \le X.
$$
We denote by $ \hat{\lambda}_n(\xi)$ the supremum  of the real numbers
$\lambda$ such that, for
any sufficiently large real number $X$, the inequalities
$$
0 < |x_0| \le X, \qquad \max_{1 \le m \le n} \,
|x_0 \xi^m - x_m| \le X^{-\lambda},
$$
have a solution in integers $x_0, \ldots, x_n$. We also set
$
\hat{w}'_n(\xi) = {1/ \hat{\lambda}_n(\xi)} .
$

In Part 2, we establish various relations linking our
six exponents.
Some of them are reformulations of known results while others are new,
as  Theorem 2.1 below. The sequel of the article is devoted
to the explicit
determination  of the values of  these six exponents in degree
$n=2$, when $\xi$ is a Sturmian continued fraction.
Since these values differ from those obtained for
algebraic numbers, we thus
obtain a quantitative  proof of the   transcendence of these
Sturmian continued fractions, which completes and precises the 
original proof in  [\ADQZ],
based on the observation that $w^*_2(\xi) > 2$ for those continued 
fractions $\xi$.
Let us indicate here that the
Sturmian words comprise  Fibonacci word and that  we thus extend
some of the results contained in
  [\RoyA, \RoyB, \ArRo] ; the above point $(2)$   being an example.
    The  formulas obtained are
    stated in  Part 3, our main  result being
  Theorem 3.1.
    Its proof, which rests on a study of multiplicative
recurrences (Lemma 4.1)  and on combinatorial properties
of infinite Sturmian words (Lemma 5.3),
is detailled in Parts 6 and 7. It extends in  a
Sturmian frame the ideas
developped in [\RoyB]. One of the key result is Lemma 6.1,
which
provides a precise estimate  of the  height of some
  quadratic numbers defined by purely periodic   continued fractions. Finally,
Part 8 is devoted to a  discussion on the  spectra
of the various exponents of approximation that we have just  defined.

In a forthcoming work, we shall combine some of the  ideas of
the present article with the method developped by
Roy [\RoyB, \RoyC, \RoyD] in order to construct
further examples of real numbers which are  badly approximable
by  algebraic integers with degree  $\le 3$.
We shall further establish a link with questions
of inhomogeneous Diophantine Approximation.

\section{2. Properties of the exponents of approximation.}

This Section is devoted to an overview of general
results on the  six functions $w_n,w_n^*,\lambda_n,
\hat{w}_n,\hat{w}_n^*$ and
$\hat{\lambda}_n$. Their values are connected by various
numerical inequalities.
Notice that   the exponents  `hat'  are uniformly
bounded in term of   $n$.  The exact determination of the upper bounds
is an important problem towards the Wirsing Conjecture
or connected topics, such as the approximation by algebraic integers; 
it is solved only for
$n=1$ and $n=2$. Let us begin with some easy properties:

\proclaim Proposition 2.1. For any integer $n \ge 1$ and any real number
  $\xi$ which is not algebraic of degree
$\le n$, we have
$$
n \le \hat{w}_n (\xi) \le w_n (\xi)
, \qquad
{1\over n} \le\hat{\lambda}_n(\xi) \le \min\{1,\lambda_n(\xi)\}
$$
and
$$
    1
\le \hat{w}^*_n (\xi) \le \min\{ w_n^* (\xi) ,\hat{w}_n (\xi)\}
\le \max\{ w_n^* (\xi) ,\hat{w}_n (\xi)\}\le w_n (\xi)
. \eqno (3)
$$

\proof The upper bound $w_n^* (\xi) \le w_n (\xi)$ is Proposition 3.2
from [\BuLiv]. The same argumentation gives also that $\hat{w}_n^*
(\xi) \le \hat{w}_n (\xi)$.
The upper bounds $\hat{w}_n (\xi) \le w_n
(\xi)$,
$\hat{w}^*_n (\xi) \le w_n^* (\xi)$ and $\hat{\lambda}_n (\xi) \le
\lambda_n (\xi)$
are easy, while the lower bounds $ \hat{w}_n (\xi)\ge n$ and
$\hat{\lambda}_n (\xi) \ge 1/n$
follow from Dirichlet's box principle (or, equivalently, Minkowski's theorem).
The lower bounds $\hat{\lambda}_1 (\xi) \ge 1$,
$\hat{w}_1 (\xi) \ge 1$ and
$\hat{w}^*_1 (\xi) \ge 1$ are direct consequences of the usual 
Dirichlet theorem
([\BuLiv], Theorem 1.1) which asserts that for any
  irrational number $\xi$ and any  $Q\ge 1$, there exists a
rational  $p/q$, with $1 \le q \le Q$ and  $0 < |\xi - p/q| <1/(qQ)$.
Since Dirichlet's box principle cannot be improved with respect to the 
exponent of $Q$
(cf. [\DaScC]), it turns out that we have the equalities
$$
\hat{w}_1 (\xi) = \hat{w}^*_1 (\xi) = \hat{\lambda}_1 (\xi)
=1.
$$
Moreover, we obviously have
$
    \hat{w}^*_n (\xi) \ge \hat{w}^*_1 (\xi)=1
$
and
$
\hat{\lambda}_n (\xi) \le \hat{\lambda}_1 (\xi) =1.
$
\cqfd

On the other hand, when  $\xi$ is a Liouville number, we have
$\hat{w}_n^*(\xi) =1$ using Liouville's inequality 
(cf. [\BuLiv], Appendix A);
it follows that the lower bound $\hat{w}_n^*(\xi) \ge 1$ in (3) is optimal.

\medskip

Wirsing proved the inequality
$$
w_n^*(\xi) \ge {w_n (\xi) \over w_n (\xi) - n + 1},  \eqno (4)
$$
  for any integer $n \ge 1$ and any real number $\xi$ which is not
algebraic of  degree at most $n$.
We can refine (4) in the following way:

\proclaim Theorem 2.1. Let $n$ be an integer $\ge 1$.  The lower bounds
$$
\hat{w}_n^* (\xi) \ge {w_n (\xi) \over w_n (\xi) -n +1}  \eqno (5)
$$
and
$$
{w}_n^* (\xi)\ge {\hat{w}_n (\xi) \over \hat{w}_n  (\xi) -n +1}, \eqno
(6)
$$
hold for any real number $\xi$ which is not
algebraic of degree  $\le n$.

\proof We obtain the first inequality by taking again Wirsing's 
argumentation [\Wir].
Let $n \ge 2$ be an integer and let $\xi$ be
a real number which is either transcendental, or algebraic of degree $> n$.
Let $\eps > 0$ and set $w = w_n(\xi) (1 + \eps)^2$. Minkowski's theorem implies
that there exist a
constant $c$ and, for any positive real number $H$, a non zero integer polynomial
$P(X)$ of degree at most $n$ such that
$$
|P(\xi)| \le H^{-w}, \quad |P(1)|, \ldots, |P(n-1)| \le H \quad
{\rm and} \quad |P(n)| \le c H^{w-n+1}. \eqno (7)
$$
The definition of  $w_n(\xi)$ and the first condition of (7) show 
that $H(P) \gg H^{1 +
\eps}$.
Consequently, $P(X)$ has some (necessarily
real) root in the neighbourhood of each of the points $\xi$, $1, \ldots, 
n-1$. Denoting by $\al$
the closest root to $\xi$, we get
$$
|\xi - \al| \gg \ll {|P(\xi)| \over H(P)} \ll H(P)^{-1} \,
(H^{w-n+1})^{-w/(w-n+1)} \and  H(P) \ll H^{w-n+1}.
$$
Since all of this is true for any sufficiently large $H$, we get 
$\hw_n^*(\xi) \ge
w/(w-n+1)$.  Selecting now $\eps$ arbitrarily close to $0$, we obtain
(5).

In order to prove (6), we take again the argumentation of [\BuGr],
reproduced in [\BuLiv], which shows that if
  $A$ is a real number, lying  in the interval $2 < A < n+1$, and   such that
$|\xi - \al| \ge H(\al)^{-A}$ for all algebraic numbers $\al$
of degree  $\le n$ and sufficiently large height, then, for
any $\eps < (n+1 - A)/(A-2)$ and any sufficiently large $H$, there
exists an integer polynomial  $P(X)$ with height  $\le H$ satisfying
$|P(\xi)| \le H^{-n-\eps}$. This means that we have
  $\hw_n (\xi) \ge n + \eps$,
and next that  $\hw_n (\xi) \ge n + (n+1 - A)/(A-2)$. Noticing that 
$w_n^*(\xi) \ge A - 1$,
we obtain (6).  \cqfd

Since ${w}_n^* (\xi) \ge \hat{w}_n^* (\xi)$
and $w_n (\xi) \ge \hat{w}_n (\xi)$, each of the two
lower bounds  from Theorem 2.1 implies the inequality (4).
Moreover,  the example of Liouville numbers $\xi$, for which
$\hat{w}_n^* (\xi) = 1$ and $w_n(\xi) = + \infty$, shows that the
estimation (5) is best possible [\LauA].

We indicate now some transference results relying the
  rational simultaneous approximation to $\xi, \ldots, \xi^n$
and the smallness  of the linear form  with coefficients $\xi,
\ldots, \xi^n$.

\proclaim Theorem 2.2.
Let $n$ be an integer $\ge 1$. We have the estimations
$$
{n \over w_n (\xi) - n + 1} \le w'_n (\xi) \le
{(n-1) w_n (\xi) + n \over w_n (\xi)} \le n \eqno (8)
$$
and
$$
{n \over \hat{w}_n (\xi) - n + 1} \le \hat{w}'_n (\xi) \le
{(n-1) \hat{w}_n (\xi) + n \over \hat{w}_n (\xi)} \le n, \eqno (9)
$$
for any real number $\xi$ which is not
algebraic of  degree  $\le n$.

\proof The inequalities (8) follow direcly from Khintchine's 
transference principle
  (cf. Theorem B.5 from [\BuLiv]),
whose proof shows in fact that (9) is also true.
\cqfd

We know the behaviour of these six functions
for almost all real numbers $\xi$, and also their values for algebraic 
real numbers $\xi$.

\proclaim Theorem 2.3. For almost all (with respect to
Lebesgue measure) real number $\xi$ and any positive integer $n$, we have
$$
w_n(\xi) = w_n^*(\xi) = w'_n(\xi) = \hat{w}_n (\xi) = \hat{w}'_n (\xi) =
\hat{w}_n^* (\xi) = n.
$$

\proof A result due to Sprind\v zuk [\SprLiv] asserts that
  $w_n (\xi)$ is equal to $n$ for almost all $\xi$.
Thank to  Proposition 2.1, to
(4) and to Theorem 2.2, we obtain that
$w_n^* (\xi) = \hat{w}_n^* (\xi) = w'_n (\xi) = n$ for almost all
$\xi$.
Besides, Bugeaud [\BuJNT] proved that $\hat{w}_n (\xi) = n$
(and therefore $\hat{w}'_n (\xi) = n$, by (9)) for almost all $\xi$. \cqfd

\proclaim Theorem 2.4. Let $\xi$ be a real algebraic number
    of degree $d \ge 2 $ and let   $n$ be a positive integer. We have
$$
w_n(\xi) = w_n^*(\xi) = w'_n(\xi) = \hat{w}_n (\xi) = \hat{w}'_n (\xi) =
\hat{w}_n^* (\xi) = \min\{n, d-1\}.
$$

\proof This is a straightforward  consequence of Schmidt's Subspace Theorem.
We omit the details. \cqfd

Davenport \& Schmidt [\DaScB] have investigated
the approximation to a real number by algebraic integers, using  a dual approch
which may also be applied to Wirsing's Conjecture. Their results have 
been recently  improved by
Laurent  [\LauB].   For a positive real number $x$, we denote by 
$\lceil x \rceil $
the smallest integer  $\ge x$. The following statement translates the 
main results
of [\DaScB] and of
[\LauB] in terms of exponents:

\proclaim Theorem 2.5. For any  integer $n \ge 1$ and any transcendental
real number $\xi$, we have
$$
w_n^* (\xi) \ge \hat{w}'_n (\xi) \ge \lceil n/2  \rceil
\and  \hat{w}_n (\xi) \le 2n - 1.
$$

Notice that the estimation
$$
w_n^* (\xi) \ge \hat{w}'_n (\xi) \ge {n \over \hat{w}_n (\xi) - n +1},
$$
which follows from  (9) and from Theorem 2.5, is only slightly weaker 
than the inequality (6).

For $n=2$, recent results of  Roy [\RoyA, \RoyB]
and Arbour \& Roy [\ArRo] give us the
  extremal values of the  functions $\hat{w}_2$, $\hat{w}'_2$
and $\hat{w}_2^*$.

\proclaim Theorem 2.6. In degree $n=2$, the refined inequalities
$$
\hat{w}'_2 (\xi) \ge {1 + \sqrt{5}\over 2}, \qquad
\hat{w}_2 (\xi) \le {3 + \sqrt{5}\over 2}
\and
\hat{w}_2^* (\xi) \le {3 + \sqrt{5}\over 2}
$$
hold for any real number $\xi$ which is neither
rational, nor quadratic. Moreover these inequalties
  are best possible.

In Theorem 2.6, the equalities are reached for
instance with the number
$$
\xi = [0; 1, 2, 1, 1, 2, 1, 2, \ldots],
$$
whose sequence of partial quotients 
reproduces the letters of the Fibonacci word, up to the initial term.
This word is an  example of
Sturmian word which will be considered in  Part 3. It corresponds to
the particular case of a constant  sequence
  $(s_k)_{k\ge 1} $  equal to one.

\section{3. Sturmian continued fractions.}

Our goal is to construct other examples of real numbers $\xi$
(necessarily transcendental)
for which one can  compute the  values  at $\xi$ of the
six functions  $w_2, \ldots, \hat{w}_2^*$, all these values
being distinct from $2$. We show that  typical Sturmian continued  fractions
share this property.

Let $(s_k)_{k\ge 1}$ be a sequence of integers $\ge 1$ and let
$\{a, b\}$ be an  alphabet with two letters. We denote by $\{a,b\}^*$
the mono\"{\i}d of  finite words on this alphabet and we define
inductively a sequence of words $(m_k)_{k \ge 0}$ belonging to $\{a,b\}^*$ by
the formulas
$$
m_0 = b, \quad m_1 = b^{s_1 - 1} a  \and
m_{k+1} = m_k^{s_{k+1}} \, m_{k-1} \qquad (k \ge 1).
$$
This sequence converges (for the product of the discrete topologies) 
in the completion
$\{a,b\}^*\cup\{a,b\}^\bN$ to the infinite word
$$
m_{\varphi} := \lim_{k \rightarrow \infty} \,
m_k = b^{s_1 - 1} a \ldots,
$$
which is usually called the {\it Sturmian characteritic word
of angle} (or of  `slope')
$$
\varphi := [0; s_1, s_2, s_3, \ldots]
$$
constructed on the  alphabet $\{a, b\}$.

The characters $a$ and
$b$ will indicate   either the letters $\{a,b\}$ of the alphabet, or
distinct positive integers, according to the context. Denote by
$\xi_{\varphi}=[0;m_\varphi]$ the real number whose partial quotients
are succesively
  $0$, and the letters of the infinite word $m_{\varphi}$, and set
$$
\sigma_{\varphi} = \liminf_{k} \, [0; s_k,s_{k-1}, \dots, s_1]  = {1
\over
\limsup_{k} \, [s_k; s_{k-1}, \dots, s_1]}.
$$
In order  to shorten the notations, we write $\sigma$
instead of
$\sigma_{\varphi}$.
\footnote{(*)}{It is readily seen that the quantity $\sigma$
does not depend upon the initial values of the sequence $(s_k)$.
Thus, in  the definition of $\sigma$,
the use  of the reverse continued fraction
$[s_k,\dots , s_1]$, ending exactly by $s_1,s_2,\dots$, is  purely 
conventional.}

The real number $\xi$ defined after  Theorem 2.6 is equal
to  the number $\xi_{(\sqrt{5} - 1)/2}$  written on the alphabet $\{1, 2\}$.
The corresponding Fibonacci word $m_{(\sqrt{5} - 1)/2} = 12112121 \ldots$
satisfies two important properties. Firstly, it   contains many
initial factors of the shape $u u u'$,
where $u$ is a word  and $u'$ is some  prefix of $u$. Since the 
quadratic real numbers
are those whose continued fraction expansion is ultimately periodic, it
follows that  $\xi_{(\sqrt{5} - 1)/2}$ has a lot of very good
quadratic approximations.
Secondly, the word $m_{(\sqrt{5} - 1)/2}$
contains many palindromic prefixes.
As was observed by Roy [\RoyA],
this last property enables us to construct rational simultaneous approximations
to $\xi_{(\sqrt{5} - 1)/2}$ and its square. It turns out that similar 
properties
also hold true for more general Sturmian words,  as we shall see in Part 5.

Here is our main result.

\proclaim
Theorem 3.1. Let $a$ and $b$ be distinct positive integers. We have
$$
\displaylines{
\lambda_2 (\xi_{\varphi}) = {1 }, \qquad
w_2(\xi_{\varphi}) = w_2^*(\xi_{\varphi}) = 1 + {2 \over \sigma},
\cr
\hat{\lambda}_2(\xi_{\varphi})= {1 +\sigma \over 2 +\sigma}= {1\over 
2} + {\sigma\over 2(2+\sigma)}
\and
\hat{w}_2(\xi_{\varphi})=\hat{w}_2^*(\xi_{\varphi})= 2 +\sigma.
\cr}
$$

It follows immediately from  Theorem 2.4 that $\xi_{\varphi}$
is a  transcendental number.  This result was originally established 
in [\ADQZ].

The proof of  Theorem 3.1 requests several steps.
We first have to control the growth of the
denominators of some convergents of $\xi_{\varphi}$ (Lemma 4.1).
Next we establish
combinatorial properties of the word $m_{\varphi}$
(Lemma 5.3). Part 6 is devoted to the computation of the exponents $w_2$,
$w_2^*$, $\hat{w}_2$ et $\hat{w}_2^*$, while
$\lambda_2$ and $\hat{\lambda}_2$ are treated in Part 7.

We briefly indicate two immediate consequences of Theorem 3.1.

\proclaim
Corollary 3.1. For any positive integer $d$, there exists a real number
$\xi$ satisfying
$$
\displaylines{
w_2(\xi) = w_2^*(\xi) = \sqrt{d^2 + 4} + d + 1,
\cr
\hat{\lambda}_2(\xi)= {\sqrt{d^2 + 4} + 2 - d \over \sqrt{d^2 + 4} +
4 - d} \and
\hat{w}_2(\xi)=\hat{w}_2^*(\xi)= 2 + {\sqrt{d^2 + 4} - d \over 2}.
\cr}
$$

\proof Apply the theorem to the constant sequence $s_k$ equal to
$d$, for which we have $\sigma = 2/(d + \sqrt{d^2 + 4})$.
\cqfd

  Corollary 3.1 gives us an explicit sequence of numbers
$>2$ (resp. $>1/2$) belonging to the spectrum of the function $\hat{w}_2$
(resp. $\hat{\lambda}_2$). More results on these spectra are obtained
in Part 8.

\proclaim
Corollary 3.2. There exists a continuous set of real numbers  $\xi$ satisfying
$$
\hat{w}'_2(\xi) < 2, \quad \hat{w}_2(\xi) > 2 \and
\hat{w}^*_2(\xi) > 2.
$$

It proceeds from Theorem 2.3 that the set of real numbers $\xi$ for which
$\hat{w}_2^* (\xi)$ is larger than $2$ has a
  Lebesgue measure equal to zero. This set  is not countable, and
it is an interesting problem to compute its
Hausdorff dimension.

\medskip

To conclude this Section, let us point out that
the  quantity $\sigma^{-1}$ appears further in  measures
of  irrationality  of some series associated to
  Sturmian sequences, as was indicated to us by Boris Adamczewski.
Denote by $(\ell_{j})_{j \ge 1}$ the sequence formed by the letters
of the word $m_{\varphi}$ written on the alphabet
$\{a, b\}$. If we choose the values $a=1$ and $b=0$, the irrationality
measure $w_1 (\beta)$ of the number
$$
\beta = \sum_{j \ge 1} \, {\ell_j \over 2^j}
$$
is exactly equal to $\limsup_{k} \, [s_k;s_{k-1}, \dots, s_1]$.
This result has been established by  Davison [\Dav] for
$\varphi= (\sqrt{5}-1)/2$ and by Adams \& Davison [\AdDa]
in the general  case (they do not state it, but the
result is implicitely contained in their work, as well as in an article
of Queff\'elec [\Que]).

\section
{4. Multiplicative recursions.}

Let $(X_n)_{ n\ge 0}$ be a sequence of positive
  real numbers and let $(s_n)_{n\ge 1}$ be a sequence of positive integers.
We suppose that the sequence $(X_n)_{ n\ge 0}$ tends to infinity
and that there   exists  $c\ge 1$ such that the inequalities
$$
c^{-s_{n+1}}  X_n^{s_{n+1}}X_{n-1} \le X_{n+1} \le c^{s_{n+1}}
X_n^{s_{n+1}}X_{n-1},
\eqno{(10)}
$$
hold  for any $n\ge 1$.
The aim of this Section
is to estimate the growth of the  sequence $(X_n)_{ n\ge 0}$.

\proclaim
Lemma 4.1. For any $n\ge 0$, define $\eta_n$ by
$X_{n+1} = X_n^{\eta_n}$. Then, we have
$$
\eta_n = [s_{n+1}; s_n, \dots , s_1](1 + \io(1))
$$
when $n$ tends to infinity.

\proof Using the right upper bound of $(10)$, we obtain successively
$$
X_n \le  c^{ s_n}X_{n-1}^{s_n}X_{n-2}
\le  c^{s_n s_{n-1} +  s_n }X_{n-2}^{s_n s_{n-1} +1}X_{n-3}^{s_n}
\le \cdots
$$
Here occurs the sequence of integers $q_k =q_{k,n}$
defined for $k= -1, \dots , n$ by
$$
q_{-1} =0,\quad  q_0 =1,\quad q_1 = s_n, \quad q_2 =s_ns_{n-1}+1, \dots
$$
and satisfying  the binary recurrent relation
$
       q_{k+1}= s_{n-k} q_k + q_{k-1}
$
       for $0\le k\le n-1$. Let us check by induction on
$k=0, \dots , n-1$ the upper bound
$$
X_n \le c^{  q_k +  q_{k-1}-1 } X_{n-k}^{q_k} X_{n-k-1}^{q_{k-1}}.
$$
The estimation is obvious for $k=0$. Suppose that it holds up to the rank $k$ ;
then  (10) gives
$$
X_n \le  c^{   q_k + q_{k-1}-1 }
\left(c^{ s_{n-k}}X_{n-k-1}^{s_{n-k}}X_{n-k-2}\right)^{q_k}
X_{n-k-1}^{q_{k-1}}
=c^{ q_{k+1}+  q_k -1} X_{n-k-1}^{q_{k+1}} X_{n-k-2}^{q_k},
$$
from which follows the upper bound at the rank $k+1$. Using the left lower
bound from (10), we obtain similarly
$$
X_n \ge c^{ -q_k -  q_{k-1} + 1} X_{n-k}^{q_k} X_{n-k-1}^{q_{k-1}}.
$$
    We shall compare $X_{n+1}$ with  $X_n$ using both estimations.
To that purpose, let us first prove the
  formula
$$
[s_{n+1}; \dots , s_{n+1-k}] = {q_{k+1,n+1}\over q_{k,n}},
$$
that is given by the theory of  continued fractions. The  recurrent relations
$$
q_{j+1,n+1}= s_{n+1-j} q_{j,n+1} + q_{j-1,n+1}
\and
q_{j,n}= s_{n+1-j} q_{j-1,n} + q_{j-2,n},
$$
are satisfied for  $1\le j \le  n$ and imply the matrix relation
$$
\left(\matrix{q_{k+1,n+1}&q_{k,n+1}\cr q_{k,n}&q_{k-1,n}\cr}\right)
=
\left(\matrix{s_{n+1}& 1 \cr 1 & 0 \cr}\right) \cdots
\left(\matrix{s_{n+1-k}& 1 \cr 1 & 0 \cr}\right),
$$
from which follows the expected formula.

We have  established the  estimations
$$
X_n \ge (c^{-1}X_{k})^{q_{n-k,n}}(c^{-1}X_{k-1})^{q_{n-k-1,n}}
\and
X_{n+1} \le (c X_{k})^{q_{n-k+1,n+1}}(c X_{k-1})^{q_{n-k,n+1}}
$$
for any $k=0, \dots , n$. Let $\epsilon $ be a positive real
number and let $k$
be a  sufficiently large integer such that  $c \le \min \{X_k^\epsilon, 
X_{k-1}^\epsilon\}$. Then, we get
$$
\eqalign{
X_{n+1} \le &
\left(c^{-1}X_{k}\right)^{{1+\epsilon\over
1-\epsilon}q_{n-k+1,n+1}
}\left(c^{-1}X_{k-1}\right)^{{1+
\epsilon\over 1-\epsilon} q_{n-k,n+1}}
\cr
\le &
X_n^{ {1+\epsilon\over 1-\epsilon}\max\left\{{q_{n-k+1,n+1}\over
q_{n-k,n}}, {q_{n-k,n+1}\over q_{n-k-1,n}} \right\} }.
\cr}
$$
Since
$$
{q_{n-k+1,n+1}\over q_{n-k,n}} =   [s_{n+1}; \dots
,s_{k+1}]=[s_{n+1}; \dots ,s_{1}]+\io_k(1)
$$
and
$$
{q_{n-k,n+1}\over q_{n-k-1,n}} =   [s_{n+1}; \dots
,s_{k+2}]=[s_{n+1}; \dots ,s_{1}]+\io_k(1)
$$
when $n$ tends to infinity, we obtain
$$
\eta_n \le [s_{n+1}; \dots ,s_{1}](1+\io(1)).
$$
The lower bound for $\eta_n$ is  proved in the same way. \cqfd

\section
{5. Combinatorial properties of Sturmian words.}

We take again the notations of Part 3, where
$m_{\varphi} $ denotes the
limit of the sequence of words  $(m_k)_{k\ge 0}$ defined by
$$
m_0= b, \quad m_1= b^{s_1-1}a \and m_{k+1}= m_k^{s_{k+1}}m_{k-1}
$$
for any $k\ge 1$. In this Section, we prove two  combinatorial properties
satisfied by the infinite word $m_{\varphi}$,
namely the existence of many (large)
initial powers and that of many initial palindromes.

For any word $u$ having at least two letters, we denote by
  $u'$ the word $u$ deprived of  its two last letters.
Thus, when  $k\ge 2$, we write  $m'_k$ for the
word $m_k$ deprived of  its two last letters
$f_k= ab$ if $k$ is even and $f_k=ba$ if $k$ is odd.

We start with a classical result on Sturmian words.

\proclaim Lemma 5.1. For any integer $k \ge 3$, we have
$m_k m'_{k-1} = m_{k-1} m'_k$ and the word
$m_{k+2}$ begins by $m_k^{1 + s_{k+1}} \, m'_{k-1} \, f_k$.

\proof This is Proposition 1 from [\ADQZ]. The second assertion
is
proved in this article, but an assertion slightly weaker is stated. \cqfd

\proclaim Lemma 5.2. For any integer $k \ge 3$, the
two words $m_{\varphi}$ and $m_k m_k\dots$ have exactly for common initial
part the word $m_k^{1+s_{k+1}}m'_{k-1}$.

\proof It follows from Lemma 5.1 that $m_k^{1 + s_{k+1}} \, m'_{k-1} \, f_k$
is a prefix of the word $m_{\varphi}$, while the periodic word
$m_k m_k \dots$ begins with
$m_k^{1+s_{k+1}} m_{k-1} = m_k^{1+s_{k+1}} m'_{k-1} f_{k-1}$.
Since $(f_{k-1}, f_k) = (ab, ba)$ or $(ba, ab)$, the lemma is established.
  \cqfd

To any word $m$ we associate the sequence, ordered by
inclusion, of prefixes of this word which are
  palindromes. The following result describes some subsequence
extracted from this sequence in the case of the Sturmian word
$m_{\varphi}$.
It should be interesting to know whether we  obtain  in this way all 
palindromic prefixes
of the word $m_\varphi$ up to a finite number of exceptions,
which seems to be experimentally the case.
  By convention, the empty word is viewed as a palindrome.

\proclaim
Lemma 5.3. The sequence of words
$$
(m_k^t m_{k-1})', \quad \hbox{\rm o\`u}\quad k \ge 1 \and 1\le t \le
s_{k+1} ,
$$
and where the  indices $(k,t)$ are lexicographically ordered,
makes up an increasing sequence of prefixes of the
word
$m_{\varphi}$ which are palindromes.

\proof
It will be convenient to state the
  relation of commutation $m_k m'_{k-1} = m_{k-1} m'_k$
of  Lemma 5.1 in the equivalent form
$$
(m_k m_{k-1})' =( m_{k-1} m_k)'
$$
which has the advantage to be valid for $k=1$ and $k=2$.

We show by induction on the rank of the couples $(k,t)$
lexicographically ordered that the  prefixes
$(m_k^t m_{k-1})'$ are palindromes.
When $k=1$ and $t\ge 1$, or $k=2$ and $t=1,2$, the words
$$
\eqalign{
(m_1^t m_{0})'= & (b^{s_1-1}a)^{t-1}b^{s_1-1} ,
\cr
(m_2m_1)'= & (b^{s_1-1}a)^{s_2}b^{s_1-1} ,
\cr
(m_2^2m_1)' = & 
b^{s_1-1}(ab^{s_1-1})^{s_2-1}ab^{s_1}a(b^{s_1-1}a)^{s_2-1}b^{s_1-1},
\cr}
$$
are clearly palindromes.
Suppose now that the recurrence hypothesis is
satisfied up to the  rank $(k,t)$ with $(k,t) \ge (2,2)$.
In order to show that
$$
(m_k^{t+1}m_{k-1})'=m_k^{t}(m_k m_{k-1})'
$$
is a palindrome, let us write
$$
\eqalign{
m_k^{t}(m_k m_{k-1})' & = (m_{k-1}^{s_k}m_{k-2})'f_{k}m_k^{t-1}
m_{k-1}^{s_k}(m_{k-2} m_{k-1})' \cr
& = (m_{k-1}^{s_k}m_{k-2})'f_{k} (m_k^{t-1}
m_{k-1})' f_{k-1} m_{k-1}^{s_k - 1}( m_{k-1} m_{k-2})' \cr
& = (m_{k-1}^{s_k}m_{k-2})'f_{k}( m_k^{t-1}
m_{k-1} )' f_{k-1} (m_{k-1}^{s_k}m_{k-2})' . \cr}
$$
Thank to the recurrence hypothesis, the words
$(m_{k-1}^{s_k}m_{k-2})'$  and $ (m_k^{t-1} m_{k-1})'$ are
palindromes,
and  $f_{k-1}$ is the mirror of $f_{k}$. We deduce that the word
$(m_k^{t+1}m_{k-1})'$ is also a palindrome.

When $t=s_{k+1}$, the successor of $ (m_k^tm_{k-1})'=m'_{k+1}$
is $m_{k+1}m'_k= m'_{k+1}f_{k+1}m'_k$. The mirror of this last
word is then
$$
m'_kf_km'_{k+1}=m_k m'_{k+1}= m_{k+1} m'_k,
$$
since the words  $m'_{k+1}$ are $m'_k$ are
palindromes, again by the recurrence hypothesis.
\cqfd

\section
{6. The exponents $w_2$, $w_2^*$, $\hat{w}_2$ and $\hat{w}_2^*$.}

We keep the  notations of Part 5 and denote by
  $\alpha_k$ the quadratic number whose continued
fraction expansion is purely
  periodic with  period $m_k$ :
$$
\alpha_k = [0; m_k, m_k,\dots].
$$
Let $M_k$ and $M'_k$ be
the  $2\times 2$ matrices associated
to the words $m_k$ and  $m'_k$ by the rule
$$
M= \prod_{i=1}^q\left(\matrix{\ell_i&1\cr 1&0\cr}\right)
$$
if the word  $m$ is written  $\ell_1\dots \ell_q$ with letters $\ell_i\in
\{a,b\}$. Recall that $a$ and
$b$ indicate at the same time letters and distinct positive integers.
Any multiplicative identity  on words translates into the same identity between
matrices $M$. In particular,
the recurrence relation $m_{k+1} = m_k^{s_{k+1}} m_{k-1}$ implies
that
$$
M_{k+1} = M_k^{s_{k+1}}M_{k-1}.
$$

If $M$ is any matrix with integer coefficients, we denote by $H(M)$ its
height, that is, the maximum of the absolute values
of its coefficients.  Let $M$ and $M'$
be $2\times 2$ matrices with positive coefficients. Then,
the height of their product
  $MM'$
clearly satisfies the relation  of quasi-multiplicativity
$$
H(M)H(M') \le H(MM') \le 2H(M)H(M'). \eqno{(11)}
$$
Let
$$
X_k =H(M_k)
$$
denote the height of the matrix $M_k$. By iterated application of
the inequalities (11), we  see that the hypothesis  (10) from Part 4
  is   satisfied  here with $c=2$.
Moreover, the sequence $X_k$ tends to infinity since the
coefficients
of $M_k$ are numerators and denominators of distinct convergents
of $\xi_\varphi$. In fact, it is easily seen that the sequence
  $X_k$ increases at least as fast as a double exponential in
  $k$, but we shall not need this point. Consequently, Lemma 4.1
can be applied to the sequence $X_k$.

In  the sequel, we write  $F\gg G$ to signify that there exists a
positive real number $\kappa$ depending only on $a , b $  and on
the sequence $s_k$, such that we have  $F(*)\ge \kappa \, G(*)$ for all
the values of the data $*$ considered in some statement.
The notation $F\gg\ll G$ signifies naturally that we
simultaneously have $F\gg G$ and $F\ll G$.

The following lemma shows that the height of the quadratic number
$\alpha_k$ has the same magnitude  as $X_k$.

\proclaim
Lemma 6.1. We have $H(\alpha_k) \gg \ll X_k$. Moreover,
the leading coefficient of the  minimal polynomial $P_{\alpha_k}$ of
$\alpha_k$ over $\bZ$ is  $\gg X_k$ and
the  conjugate  $\alpha'_k$ of  $\alpha_k$ satisfies
$|\alpha'_k - \alpha_k| \gg 1$.

\proof
The number  $1/\alpha_k$ is a fixed point of the homography 
associated to the matrix
$M_k$. The upper bound  $H(\alpha_k) \ll X_k$
is then obvious. The lower bound $H(\alpha_k) \gg X_k$
needs more elaborate arguments. Notice that
$$
M_k= M'_k F_k
$$
with
$$
F_k = \left(\matrix{a &1\cr 1&0\cr}\right) \left(\matrix{b &1\cr
1&0\cr}\right) =
\left(\matrix{ab+1 &a\cr b& 1\cr}\right)
\quad \hbox{\rm or }\quad
F_k = \left(\matrix{b &1\cr 1&0\cr}\right) \left(\matrix{a &1\cr 1
&0\cr}\right)
=\left(\matrix{ab+1 &b\cr a& 1\cr}\right),
$$
according to the parity of  $k$.
    Set $F_k=\left(\matrix{f_0 &f_1\cr f_2&f_3\cr}\right)$
    and observe that the  matrix $M'_k=\left(\matrix{x_{0}
&x_{1}\cr x_{1}&x_{2}\cr}\right)$ is  symmetrical,
since  $m'_k$ is a palindrome. The quadratic number $\alpha_k$
is then a root of the trinomial
$$
(f_1x_{0}+f_3x_{1})T^2 +(f_0x_{0}+(f_2-f_1)x_{1}-f_3x_{2})T
-f_0x_{1}-f_2x_{2}.
$$
The coefficients of the three linear forms in
$x_{0},x_{1},x_{2}$ defined by the three  coefficients
of the above  trinomial  form  a matrix whose determinant
$$
(f_2-f_1)(f_0 f_3 -f_1 f_2)=\pm(a-b)
$$
is non zero, because $a$ and $b$ are distinct.
Since the integers $x_0,x_1,x_2$ are relatively prime, the gcd
of the coefficients of this trinomial is  bounded.  It follows that
the leading coefficient of the polynomial $P_{\alpha_k}$ is obtained
by dividing $f_1x_0+ f_3x_1$ (which is $\gg X_k$) by this gcd and that
$$
H(\alpha_k)\gg H(M'_k)\gg H(M_k) = X_k.
$$

Notice finally that the discriminant of the trinomial is
$$
\ge (f_1 x_0 + f_3x_1)(f_0 x_1 + f_3 x_1) \gg X_k^2 ,
$$
since the  coefficients of the
matrices
$M'_k$ and $F_k$ are all $\ge 0$. It follows that
$|\alpha'_k - \alpha_k| \gg  1$. An other way to prove this
result rests on the observation  that the continued fraction expansion of
$1/ \alpha_k$ is purely periodic and therefore, by a classical result
of Galois (cf. for example [\RoSz], Theorem 3, page 45), its conjugate
belongs to the interval $[0, 1]$.  \cqfd

\bigskip

   We  deduce now from Lemmas 5.2 and 6.1 the exact  values of
  $w_2 (\xi_{\varphi})$ and of $w_2^*(\xi_{\varphi})$.

\proclaim Corollary 6.1.
We have  $w_2 (\xi_{\varphi}) = w_2^*(\xi_{\varphi}) = 1 + 2/\sigma$.

\proof It follows from Lemma 5.2 that
the words $m$ and $m_km_k\dots$ have
  $m_k^{1+s_{k+1}}m'_{k-1}$ as a common prefix.
Thus the partial quotients  of the numbers
$\xi_{\varphi}$ and $\alpha_k$ coincide first with $0$,
next with the letters of the word $m_k^{1+s_{k+1}}m'_{k-1} = m_km'_{k+1}$, and
differ just after. Let $N$ be the number of letters of the word $
m_k m'_{k+1}$, and let
$$
{p_0\over q_0}=0, {p_1\over q_1}, \dots ,{p_N\over q_N},
$$
be the $N+1$ first convergents which are common to both real numbers
$\xi_\varphi$
and $\alpha_k$. Then, the theory of continued fractions  shows us that
$$
M_k M'_{k+1} =  \left(\matrix{q_N &q_{N-1}\cr p_N & p_{N-1}\cr}\right) ,
$$
and that we have the estimate
$$
{1\over q_N^2} \ll \vert\xi_\varphi -\alpha_k\vert \le {1\over q_N^2},
$$
the lower bound coming from the fact that the partial quotients of the numbers
$\xi_\varphi$ and $\alpha_k$ are bounded, since they are equal to   $a$ or to
$b$.
    Thank to Lemmas 4.1 and 6.1,
we thus obtain the inequalities
$$
\vert \xi_{\varphi}-\alpha_k\vert \gg\ll H(M_kM'_{k+1})^{-2}\gg\ll
(X_kX_{k+1})^{-2}
\gg\ll H(\alpha_k)^{-2 -2 \eta_k},
\eqno (12)
$$
from which follows the lower bound
$w_2^*(\xi_{\varphi}) \ge  1 + 2/\sigma$.
In particular, we have $w_2^*(\xi_{\varphi}) = + \infty$ whenever the sequence
$(s_k)_{k \ge 1}$ is unbounded.
If this sequence is bounded, we proceed in the following way
to get an upper bound for $w_2^*(\xi_{\varphi})$.

Let $\alpha$ be a real algebraic number of degree $\le 2$.
Let  $k$ be the integer defined by the inequalities
$$
H(\alpha_{k-1}) \le H(\alpha) < H(\alpha_k).
$$
Liouville's inequality (cf. [\BuLiv], Appendix A) gives the lower
bound
$$
|\alpha_{k} - \alpha| \gg  H (\alpha)^{-2} \,
H(\alpha_k)^{-2} \gg H(\alpha_k)^{-4},
$$
while (12) implies that
$$
|\xi_{\varphi} - \alpha_k| \ll H(\alpha_k)^{-2 - 2 \eta_k}
\ll H(\alpha_k)^{-4-\epsilon}
$$
for some positive real number $\epsilon$
depending on the maximal value of the $s_k$.
Then, we deduce from the triangle inequality that
$$
|\xi_{\varphi} - \alpha| \ge |\alpha_k - \alpha| - |\xi_{\varphi} -
\alpha_k|
\gg H(\alpha)^{-2} \, H(\alpha_k)^{-2}
\gg H(\alpha)^{-2 - 2 \eta_{k-1}},
$$
which implies that $w_2^* (\xi_{\varphi}) \le 1 + 2/\sigma$.

In order to show that
$w_2 (\xi_{\varphi}) = 1 + 2/\sigma$,
one proceeds in a similar way, evaluating any integer polynomial
$P(X)$  of degree $\le 2$ at the point $\alpha_k$, where the index $k$ 
is defined by
$H(\alpha_{k-1}) \le H(P) < H(\alpha_k)$.  The details are left to 
the reader.\cqfd

\bigskip

The quadratic approximations  furnished by Lemma 6.1 are sufficiently
dense to determine exactly the values of  the exponents
$\hat{w}_2 (\xi_{\varphi})$ and $\hat{w}_2^*(\xi_{\varphi})$.

\proclaim Corollary 6.2.
We have $\hat{w}_2 (\xi_{\varphi}) =
\hat{w}_2^*(\xi_{\varphi}) = 2 + \sigma$.

\proof
Lemma 4.1, together with (12) rewritten  in the form
$$
\vert \xi_{\varphi}-\alpha_k\vert
\gg\ll (X_kX_{k+1})^{-2}
\gg\ll H(\alpha_k)^{-1}H(\alpha_{k+1})^{ -\left(2 + {1\over
\eta_k}\right)} ,
\eqno{(13)}
$$
implies  the lower bound
$$
\hat{w}_2^*(\xi_{\varphi})\ge
2+\liminf_{k}([0;s_k,s_{k-1}, \dots, s_1]) = 2+\sigma .
$$
Suppose now that
$$
\hat{w}_2 (\xi_{\varphi}) > 2+\sigma .
$$
We shall reach a  contradiction starting with an  integer $k$  and
real numbers $\epsilon , w$, such that
$$
\hat{w}_2 (\xi_{\varphi}) > w \ge 2 +{1\over \eta_k} +\epsilon.
$$
This will imply  $\hat{w}_2 (\xi_{\varphi}) \le 2+\sigma$
and therefore
$$
\hat{w}_2^*(\xi_{\varphi}) = \hat{w}_2 (\xi_{\varphi}) = 2+\sigma.
$$
By assumption, we may suppose $k$ arbitrarily large with $w$ and
$\epsilon$ fixed.
Set $X= \rho X_{k+1}$ for some small constant $\rho $ (independent
of $k$). Thus, there exists a non zero integer polynomial $P(X)$ of
degree $\le 2$ such that
$$
H(P) \le X \and \vert P(\xi_{\varphi}) \vert \le X^{-w}.  \eqno (14)
$$
Let us prove that $P(\alpha_k)= 0$.
Otherwise, Liouville's inequality, together with the triangle inequality
and Rolle's Lemma
applied to $|P(\alpha_k) - P(\xi_{\varphi})|$, gives us the inequalities
$$
\eqalign{
\rho^{-1} X_{k+1}^{-1} X_{k}^{-2} \ll H(P)^{-1}H(\alpha_k)^{-2} & \ll
\vert
P(\alpha_k) \vert \cr
& \ll \max\Big\{H(P) (X_kX_{k+1})^{-2}, (\rho X_{k+1})^{-w}\Big\}.\cr}
\eqno (15)
$$
If the maximum on the right hand side of (15) is reached by the 
second term, we obtain
$$
H(P) \gg (\rho X_{k+1})^wH(\alpha_k)^{-2}
\gg \rho^wX_{k+1}^{2+{1\over \eta_k}+\epsilon} X_k^{-2}
= \rho^{w} X_{k+1}^{2 - {1\over \eta_k}+\epsilon}.
$$
Bounding roughly from below  $\eta_k \ge 1 - o(1)$ by Lemma 4.1,
the above lower bound for $H(P)$ then contradicts its upper bound
$H(P)\le \rho X_{k+1}$, when $k$ is large enough.
Therefore, the maximum in  (15) is equal to
$H(P) X_k^{-2}X_{k+1}^{-2}$. Then (15) implies that $H(P)\gg 
\rho^{-1} X_{k+1}$,
in contradiction with the upper bound $H(P) \le  \rho X_{k+1}$ from (14), for
$\rho$ small enough .

Since the  minimal polynomial $P_{\alpha_k}(X)$
    satisfies (14), we deduce from  Lemma 6.1 and
from  (13) the estimation
$$
X_k^{-1} X_{k+1}^{-2} \ll X_k \vert
\xi_{\varphi} - \alpha'_k \vert \vert \xi_{\varphi} - \alpha_k \vert
\ll \vert P_{\alpha_k}(\xi_{\varphi}) \vert \ll X_{k+1}^{-w},
$$
which implies the upper bound
$$
w \le 2 +{1\over \eta_k} +\io(1).
$$
Let us recall  that we have selected $k$ arbitrarily large with
  $w \ge 2 +{1\over \eta_k} +\epsilon$.
The final contradiction follows. \cqfd

\section
{7. Construction of simultaneous rational approximations to
$\xi_{\varphi} , \xi_{\varphi}^2$.}

As for the Fibonacci continued fraction $\xi_{(\sqrt{5}-1)/2}$
considered by  Roy [\RoyA], the ordered sequence of palindromic prefixes of the
  mot $m_{\varphi}$ provides a sequence of simultaneous rational approximations
to $\xi_{\varphi} ,
\xi_{\varphi}^2$. Let us recall the principle of the construction.

If the  prefix $m_{\varphi\vert N}$ of the first  $N$
letters of the word $m_\varphi$ is a palindrome, then the matrix
$$
M_{\vert N} = \left(\matrix{x_0 &x_1\cr x_1& x_2\cr}\right)
$$
associated to the word  $m_{\varphi\vert N}$
is  symmetrical, and we have
$$
x_0= q_N,\quad x_1 =p_{N}=q_{N-1},\quad x_2=p_{N-1},
$$
where
$
{p_0/ q_0}=0, {p_1/ q_1},\dots
$
denotes the sequence of convergents of $\xi_{\varphi}$.
Thus,
$$
{x_2\over x_0}= {p_N\over q_N} \cdot {p_{N-1}\over q_{N-1}}
$$
is a good rational approximation of  $\xi_\varphi^2$.
Setting
$$
\ux = (x_0, x_1, x_2), \quad  H(\ux) = \max\{x_0, x_1, x_2 \}\and
L(\ux) = \max\{\vert x_0\xi_{\varphi} -x_1\vert,
\vert x_0\xi_{\varphi}^2 -x_2\vert\},
$$
it follows that
$$
L(\ux)
\gg\ll H(\ux)^{-1}. \eqno (16)
$$
With some abuse of notation, we shall identify without ambiguity
the integer triple $\ux=(x_0,x_1,x_2)$
with the symmetrical matrix $M_{\vert N}$.
Thus, we have
$$
\det(\ux) =x_0x_2-x_1^2=\pm 1 ,
$$
since the matrix  $\ux$ is defined by a product of  matrices of
determinant $-1$.
It is also convenient to view $\ux$ as the projective coordinates of a point
in  $\bP^2$, and to investigate
  alignments of points in $\bP^2$. Three points $\ux,\uy,\uz$
are located on the same  line in $\bP^2$ if $\det (\ux, \uy , \uz) =0$.

Let us uniquely write each integer $\ell \ge 1$ in the
form
$$
\ell = \ell_k +t \quad \hbox{\rm where} \quad \ell_k = s_2 + \cdots +
s_k,
\quad (\ell_1 =0)
$$
with $k\ge 1$ and $ 1\le t \le s_{k+1}$. We denote by $\ux_{\ell}$ the
triple associated to the  palindrome
$(m_k^tm_{k-1})'$ of index $(k,t)$ considered in Lemma 5.3.
The following result describes the configuration in $\bP^2$ of the sequence of
  points $\ux_\ell$.

\goodbreak

\medskip
\sl
\noindent
{\bf Lemma 7.1.} Let $k$ be an  integer $\ge 3$.

\noindent
(i) \quad The point $\ux_{\ell_{k-1}}$ and the  $s_{k+1}+1$
consecutive points
$$
\ux_{\ell_k +1},\dots , \ux_{\ell_{k+1}+1}
$$
are located on the same projective line.

\noindent
(ii) \quad We have the formula
$$
\quad \quad \det (\ux_{\ell_{k}}, \ux_{\ell_k +1},\ux_{\ell_k +2}) =
\pm (b-a).
$$
Therefore, the  three consecutive points  $ \ux_{\ell_k},\ux_{\ell_k
+1 },\ux_{\ell_k +2}$
are linearly independent.
\rm
\bigskip

\proof
Notice that the point $\ux_{\ell_{k-1}}$
is associated to the  palindrome $m'_{k-1}$.
When $k\ge 3$ and $\ell =\ell_k+t$, with $1\le t\le s_{k+1}$, the
$\ell$-th palindrome
considered in Lemma 5.3 is
equal to $m_k^tm'_{k-1}$. In terms of matrices $M_j$ and $M'_j$, we
have
$$
\ux_{\ell_{k-1}} = M'_{k-1} ,\quad  \ux_{\ell_k +t} = {M_k}^t M'_{k-1}
\eqno{(17)}
$$
for $t= 1, \dots ,s_{k+1}$, and
$$
\ux_{\ell_{k+1}+1} = M_{k+1}M'_k = M_k^{s_{k+1}}M_{k-1}M'_{k} =
M_k^{1+s_{k+1}}M'_{k-1}, \eqno{(18)}
$$
the last equality  coming from the relation  $m_{k-1}m'_{k}=
m_{k}m'_{k-1}$ (cf. Lemma 5.1) at the level of
the associated words. We observe that all these matrices are of the
form
$M_{k}^tM'_{k-1}$
with
$t=0,1, \dots , s_{k+1} +1$.
Then, the Cayley-Hamilton Theorem  gives
the linear relations
$$
\ux_{\ell_k +2} = \hbox{\rm trace}(M_{k}) \ux_{\ell_k +1 } \pm
\ux_{\ell_{k-1}}
\and
\ux_{\ell_k +t+2} = \hbox{\rm trace}(M_{k}) \ux_{\ell_k +t+1 } \pm
\ux_{\ell_{k}+t}
\eqno{(19)}
$$
linking three consecutive points of the list (i).

Concerning the assertion  (ii), we deduce by linearity from  (19)  that
$$
\det (\ux_{\ell_{k}}, \ux_{\ell_k +1},\ux_{\ell_k +2}) =
\pm\det (\ux_{\ell_{k}}, \ux_{\ell_k +2},\ux_{\ell_k +3})
=\cdots =
\pm \det (\ux_{\ell_{k}}, \ux_{\ell_{k +1}},\ux_{\ell_{k+1} +1}).
$$
Now, we have
$$
\ux_{\ell_{k}} = M'_k , \quad \ux_{\ell_{k +1}} =M'_{k+1} \and
\ux_{\ell_{k+1} +1} = M_{k+1}M'_k.
$$
Following Roy,
let us set  $J=\left(\matrix{0 &1\cr -1&0\cr}\right)$ and
use his formula (2.1) from [\RoyB]:
$$
\det(\ux ,\uy ,\uz) = \hbox{\rm trace }(J\ux J \uz J  \uy),
$$
which is valid for all
triples $\ux,\uy,\uz$, identified as before to symmetrical
matrices.
Noting that  $J\ux J= \pm \ux^{-1}$ if $\det(\ux)=\pm 1$, we obtain
$$
\eqalign{
\det (\ux_{\ell_{k}}, \ux_{\ell_{k +1}},\ux_{\ell_{k+1} +1})
= &
\pm \hbox{\rm trace }( J M'_k (M_{k+1}M'_{k})^{-1} M'_{k+1})
\cr
= &\pm \hbox{\rm trace }( JF_{k+1}^{-1}) = \pm (b-a),
\cr}
$$
which is non zero since $a$ and $b$ are distinct. \cqfd

\bigskip

Now, we are able to find the values
of the functions $\lambda_2$ and $\hat{\lambda}_2$ at the point 
$\xi_{\varphi}$.

\proclaim Proposition 7.1. We have
$$
\lambda_2 (\xi_{\varphi}) = {1 } \and
\hat{\lambda}_2(\xi_{\varphi}) = {1 +\sigma \over 2 +\sigma}.
$$

\proof
The first formula $\lambda_2 (\xi_{\varphi})=1$
follows immediately from (16), since the partial quotients of $\xi_{\varphi}$
are bounded.

The proof of the second assertion is more elaborate.
For any $\ell \ge 1$, let us set $X'_\ell = H(\ux_\ell)$ and denote by
$\eta'_\ell$ the exponent defined by the
relation $X'_{\ell+1}= {X'_l}^{\eta'_\ell}$.

We  first prove that
$$
\liminf_{\ell\ge 1}\left({1 \over \eta'_\ell}\right) = {1+\sigma\over
2+\sigma}.
$$
One deduce from  (11) and from the relations (17) and (18), the estimation
$$
X_k^tX_{k-1} \ll X'_{\ell_k+t} \ll (2X_k)^t X_{k-1}
$$
which is valid for $t=1,\dots , s_{k+1}+1$ and for any integer $k\ge 
3$. It follows that
$$
\eta'_{\ell_k +t} = {t+1 +{1\over\eta_{k-1}}\over t
+{1\over\eta_{k-1}}} \left( 1 + \cO \left( {1\over \log
X_k}\right)\right)
\le {2 +{1\over\eta_{k-1}}\over 1 +{1\over\eta_{k-1}}} \Big( 1
+o(1)\Big),
\eqno{(20)}
$$
for  $t=1, \dots ,s_{k+1}$ and $k$ tending to infinity. Then Lemma
4.1 implies that
$$
\liminf_{\ell}\left({1 \over \eta'_\ell}\right)
=\liminf_{k }\left({1+{1\over \eta_{k}}\over 2 +{1\over
\eta_{k}}}\right) ={1+\sigma\over 2+\sigma}.
$$

     The upper bound
$$
L(\ux_\ell) \ll {X'_\ell}^{-1} \eqno{(21)}
$$
is a special case of (16). Selecting in the
interval $X'_\ell \le
X < X'_{\ell+1}$ the point $ \ux_\ell$ as a simultaneous rational approximation
of $(\xi_\varphi,\xi_\varphi^2)$, we obtain immediately the lower bound
$$
\hat{\lambda}_2(\xi_{\varphi}) \ge
\liminf_{\ell\ge 1}\left({1 \over \eta'_\ell}\right)
    ={1+\sigma\over 2+\sigma}.
$$

In order to prove the equality  $\hat{\lambda}_2(\xi_{\varphi})
=(1+\sigma)/(2+\sigma)$, we  argue now by contradiction.
The idea rests on the observation that the inferior limit
of the sequence  $1/\eta'_\ell$ is reached for indices
$\ell$ of the form  $\ell_k+1$, as (20) shows. If on the contrary,
$\hat{\lambda}_2(\xi_{\varphi})$
should be $>(1+\sigma)/(2+\sigma)$, let $\lambda$ be a real number in the
interval
$$
{1+\sigma\over 2+\sigma}
=\liminf_{\kappa}\left({1 \over \eta'_{\ell_\kappa +1}}\right)
< \lambda <
\hat{\lambda}_2(\xi_{\varphi}).
$$
Thus, there exist infinitely many integers $k$ such that
$$
     \lambda
\ge {1\over \eta'_{\ell_k+1}}  . \eqno{(22)}
$$
    Furthermore, for each  sufficiently large $X$, there exists
    some non zero integer point $\ux $ such that
$$
H(\ux) \le X \and L(\ux) \le X^{-\lambda} .
$$
Choose $X= \rho X'_{\ell_k +2}$, for some positive constant $\rho$
  independent of  $k$, and observe that
$$
X'_{\ell_k} = X_k, \quad
X'_{\ell_k+1}\gg\ll X_kX_{k-1}, \quad
X'_{\ell_k+2}\gg\ll X_k^2X_{k-1}. \eqno{(23)}
$$
Then we deduce from (22) and from (23) that
$$
    L(\ux) \le X^{-\lambda} \le \rho^{-\lambda}
{X'_{\ell_k +1}}^{-1} \ll
\rho^{-\lambda} (X_kX_{k-1})^{-1}. \eqno{(24)}
$$

  Now we prove that
$$
\det(\ux , \ux_{\ell_k +1}, \ux_{\ell_k}) = \det(\ux , \ux_{\ell_k
+1}, \ux_{\ell_k+2}) = 0 ,
$$
for  $k$ sufficiently large, if the constant $\rho$ is small enough. 
Proceeding as in
  Lemma 4 of [\DaScB], and using the estimations (21),(23) and
(24), we bound
$$
\eqalign{
\vert \det(\ux , \ux_{\ell_k +1}, \ux_{\ell_k}) \vert
\ll &
X'_{\ell_k+1}L(\ux_{\ell_k})L(\ux) +X'_{\ell_k}L(\ux_{\ell_k+1})
L(\ux) + XL(\ux_{\ell_k})L(\ux_{\ell_k+1})
\cr
\ll & \rho^{-\lambda}X_k^{-1} +\rho^{-\lambda}X_k^{-1}X_{k-1}^{-2} +
\rho
\cr}
$$
and
$$
\eqalign{
\vert \det(\ux , \ux_{\ell_k +1}, \ux_{\ell_k +2}) \vert
\ll & X'_{\ell_k+2}L(\ux_{\ell_k +1}) L(\ux) +
X'_{\ell_k+1} L(\ux_{\ell_k+2}) L(\ux)
+ X L(\ux_{\ell_k+2}) L(\ux_{\ell_k+1})
\cr
\ll & \rho^{-\lambda}X_{k-1}^{-1} + \rho^{-\lambda}X_k^{-2}X_{k-1}^{-1}
+
\rho (X_{k}X_{k-1})^{-1} .
\cr}
$$
      These two determinants being integers, both necessarily vanish
if $\rho$ is small and if the index $k$ is large enough.

  The assertion (ii) of Lemma 7.1 then implies that $\ux$ is an integer multiple
of $\ux_{\ell_k +1}$,
since this last  point has relatively prime coordinates. Thus
$$
{X'_{\ell_k+1}}^{-1} \ll L(\ux_{\ell_k +1}) \le L(\ux) \le (\rho
X'_{\ell_k +2})^{-\lambda}
=\rho^{-\lambda}{X'_{\ell_k +1}}^{-\lambda \eta'_{\ell_k +1}}.
$$
We obtain
$$
\lambda \le { 1 +\io(1)\over \eta'_{\ell_k +1}},
$$
from which follows the upper bound $\lambda \le \liminf_\kappa 
(1/\eta'_{\ell_\kappa
+1})$, since the above bound holds for infinitely many $k$.
This last upper bound contradicts our choice of
  $\lambda$. \cqfd

\section
{8. Questions of spectra.}

   By abuse of language, we call
{\it spectrum} of a  function the set of values taken
by this function at  transcendental points. In view of Theorem
2.4, this  definition is by no means  restrictive when we  consider
the six exponents of approximation studied in the present article.

We start with a brief review of  known results
about the functions $w_n$ and $w_n^*$.
Thank to exact computations of   Hausdorff dimensions,
Baker \& Schmidt [\BaSc] have shown that the spectrum of the function
$w_n^*$ contains $[n, + \infty]$. We have equality if the
Wirsing Conjecture, recalled in the introduction, holds. In a long
and difficult article, Bernik [\Bern] has established that the spectrum of
the function $w_n$ is equal to the interval $[n, + \infty]$.
As far as we know, the study of the spectrum of $\lambda_n$ remains
to be performed. Nevertheless,
  it is easily seen that this spectrum is unbounded.
  Remark also that the functions $w_n$ and
$w_n^*$ differ at  some transcendental points  [\RCBak, \BuAA].

Unlike $w_n,\, w_n^*$ and  ${\lambda}_n$, the functions $\hat{w}_n$,
$\hat{w}_n^*$ and $\hat{\lambda}_n$ are bounded and their  spectrum is,
respectively, contained in the  intervals $[n, 2n - 1]$,
$[1, 2n - 1]$ and $[1/n, (\lceil n/2  \rceil)^{-1}]$,
for any degree $n \ge 1$.

Now, let us look more carefully at the spectra of the
functions $\hat{w}_2$ and $\hat{\lambda}_2$,  respectively contained in
the more restrained intervals $[2,(3+\sqrt{5})/2]$ and
$[1/2,(\sqrt{5}-1)/2]$, according to
  Theorem 2.6.  Our Theorem 3.1
   points out a lot of values in each of these quadratic spectra, and it
invites us to study the set $\cS$ of
values taken  by the quantity
$$
    \sigma = {1 \over \displaystyle
    \limsup_{k \ge 1} \, [s_k; s_{k-1}, \ldots, s_1] } ,
$$
for any bounded sequence $s_k$  of positive integers.
It is easily seen that the largest element  in $\cS$ is
$( \sqrt{5}-1)/2 \simeq 0.618...$ which is
obtained when  $1$ is the only integer appearing
infinitely many times in the sequence $s_k$.
The values immediately below have been found
by  Cassaigne [\Cassa]. They constitute a decreasing sequence of
quadratic numbers converging to
the largest accumulation point
of $\cS$. His result, formulated in Theorem 2 of
[\Cassa], is the following:

\proclaim
Theorem 8.1.
Denote by  $\psi$ the endomorphism of the mono\"\i d
    $\{1,2\}^*$ defined by
$$
\psi(1) =2 \and \psi(2) = 211.
$$
The decreasing sequence of quadratic numbers
$$
\sigma_n = [0;\psi^n(1),\psi^n(1),\dots ] \qquad (n\ge 0)
$$
converges to the real number
$$
s= [0; \lim_{n\rightarrow \infty}\psi^n(1)]=
[0;2,1,1,2,2,2,1,1,\dots] \simeq   0.38674997056...
$$
The point $s$ is the largest accumulation point of $\cS$ and
the intersection of $\cS$ with the interval $]s,+\infty]$
coincides with the set of values taken by the sequence 
$(\sigma_n)_{n\ge 0}$.

\medskip

    In  Theorem 8.1, the
    sequence of iterated words $\psi^n(1)$ converges
to the infinite word  $21122211\dots$, the fixed point of the substitution
$\psi$ acting on $\{1,2\}^\bN$.
Remark also that
  $\sigma_n = \sigma_\varphi$ for the choice of angle
$\varphi=  [0;\widetilde{\psi^n(1)},\widetilde{\psi^n(1)},\dots ]$,
where the symbol $ \,\,\widetilde{} \,\,$
stands for the mirror of a word.

Thus the three largest values of  $\hat{w}_2(\xi)$, when $\xi$
is any Sturmian continued fraction, are
$$
{3+\sqrt{5}\over 2} > 1+\sqrt{2} > {4 +\sqrt{10}\over 3} > \cdots
$$
while those of $\hat{\lambda}_2(\xi)$ are
$$
{\sqrt{5}-1\over 2} > 2 - \sqrt{2} > {-2 +\sqrt{10}\over 2} > \cdots
$$
Then, a natural  question arises: where are located these values
inside the larger  spectra
of the values $\hat{w}_2(\xi)$ and $\hat{\lambda}_2(\xi)$,
when $\xi$ now ranges  among all transcendental real numbers\thi ?
If we exclude
extremal numbers of   Roy for which we simultaneously have
$\hat{w}_2(\xi) = (3+\sqrt{5})/ 2$ and
$\hat{\lambda}_2(\xi) =(\sqrt{5}-1) / 2$, does there exist
real numbers $\xi$ such that
    $\hat{w}_2 (\xi) > 1+ \sqrt{2}$, or  $\hat{\lambda}_2(\xi) > 2
-\sqrt{2}$  ?
In  this respect, one can remark that other
  combinatorial constructions lead to real numbers $\xi$
for which $\hat{w}_2^* (\xi) > 2$ and $\hat{\lambda}_2(\xi) >
1/2$. As an example, we can consider
the number $\xi$ whose  partial quotients are given by
the  fixed point of the so-called  Tribonacci substitution defined by
$1 \to 12$, $2 \to 13$ and $3 \to 1$.

Besides, Cassaigne has shown in [\Cassa]  that the subset
  $\cS \subset \bR$ has an empty interior. It should be interesting to
know whether the same result holds for the
  spectra of the functions $\hat{w}_2$ and   $\hat{\lambda}_2$.

Let us finally point out that the study of the set $\cS$ below its
greatest accumulation point $s$ remains to be done.

\vskip 8mm

\centerline{\bf References}

\vskip 5mm

\item{[\AdDa]}
W. W. Adams and J. L. Davison,
{\it A remarkable class of continued fractions},
Proc. Amer. Math. Soc. 65 (1977), 194--198.

\item{[\ADQZ]}
J.-P. Allouche, J. L. Davison, M. Queff\'elec, and L. Q. Zamboni,
{\it Transcendence of Sturmian or morphic continued fractions},
J. Number Theory 91 (2001), 39--66.

\item{[\ArRo]}
B. Arbour and D. Roy,
{\it A Gel'fond type criterion in degree two},
Acta Arith. 11 (2004), 97-103.

\item{[\BaSc]}
A. Baker and W. M. Schmidt,
{\it Diophantine approximation and
Hausdorff dimension},
Proc. London Math. Soc. 21 (1970), 1--11.

\item{[\RCBak]}
R. C. Baker,
{\it On approximation with algebraic numbers of bounded degree},
Mathematika 23 (1976), 18--31.

\item{[\Bern]}
   V. I. Bernik,
{\it Application of the Hausdorff dimension in the theory of
Diophantine approximations}, Acta Arith. 42 (1983), 219--253 (in Russian).
English translation in Amer. Math. Soc. Transl. 140 (1988), 15--44.

\item{[\BuGr]}
     Y. Bugeaud,
{\it On the approximation by algebraic numbers
with bounded degree},
Algebraic number theory and Diophantine analysis (Graz, 1998), 47--53,
de Gruyter, Berlin, 2000.

\item{[\BuJNT]}
     Y. Bugeaud,
{\it Approximation par des nombres alg\'ebriques},
J. Number Theory 84 (2000), 15--33.

\item{[\BuAA]}
     Y. Bugeaud,
     {\it Mahler's classification of numbers compared with Koksma's},
Acta Arith. 110 (2003), 89--105.

\item{[\BuLiv]}
     Y. Bugeaud,
     Approximation by algebraic numbers,
Cambridge Tracts in Mathematics, Cambridge University Press. To appear.

\item{[\Cassa]}
J. Cassaigne,
{\it Limit values of the recurrence quotient of Sturmian sequences},
Theor. Comput. Sci. 218 (1999), 3--12.


\item{[\DaScA]}
H. Davenport and W. M. Schmidt,
{\it Approximation to real numbers
by quadratic irrationals}, Acta Arith. 13 (1967), 169--176.

\item{[\DaScB]}
     {H. Davenport and W. M. Schmidt},
{\it Approximation to real numbers by
algebraic integers}, Acta Arith. {15} (1969), 393--416.

\item{[\DaScC]}
H. Davenport and W. M. Schmidt,
{\it Dirichlet's theorem on
Diophantine approximation}, Symposia Mathematica, Vol. IV
(INDAM, Rome, 1968/69), pp. 113--132, Academic Press, London, 1970.

\item{[\Dav]}
J. L. Davison,
{\it A series and its associated continued fraction},
Proc. Amer. Math. Soc. 63 (1977), 29--32.

\item{[\Kok]}
J. F. Koksma,
{\it \"Uber die Mahlersche Klasseneinteilung der transzendenten Zahlen
und die Approximation komplexer Zahlen durch algebraische Zahlen},
Monats. Math. Phys. 48 (1939), 176--189.

\item{[\LauA]}
     M. Laurent,
{\it Some remarks on the approximation of complex numbers by
algebraic numbers}.
Proceedings of the 2nd Panhellenic Conference in Algebra and
Number Theory (Thessaloniki, 1998). Bull. Greek Math.
Soc. 42 (1999), 49--57.

\item{[\LauB]}
     M. Laurent,
{\it Simultaneous rational approximation to the successive powers
of a real number}.
Indag. Math. 11 (2003), 45--53.

\item{[\Mah]}
K. Mahler,
{\it Zur Approximation der Exponentialfunktionen und des
Logarithmus. I, II},
J. reine angew. Math. 166 (1932), 118--150.

\item{[\Que]}
M. Queff\'elec,
{\it Approximations diophantiennes des nombres sturmiens},
J. Th\'eor. Nombres Bordeaux 14 (2002), 613--628.

\item{[\RoSz]}
A. M. Rockett and P. Sz\"usz,
Continued Fractions, World Scientific, Singapore, 1992.

\item{[\RoyA]}
     D. Roy,
{\it Approximation simultan\'ee d'un nombre et son carr\'e},
C. R. Acad. Sci. Paris 336 (2003), 1--6.

\item{[\RoyB]}
     D. Roy,
{\it Approximation to real numbers by cubic algebraic numbers, I},
Proc. London Math. Soc. {88} (2004), 42--62.

\item{[\RoyC]}
     D. Roy,
{\it Approximation to real numbers by cubic algebraic numbers, II},
Ann. of Math. {158} (2003), 1081--1087.

\item{[\RoyD]}
     D. Roy,
{\it Diophantine approximation in small degree},
Preprint.

\item{[\SprLiv]}
V. G. Sprind\v zuk, Mahler's problem in metric number theory,
Amer. Math. Soc., Transl. math. monogr. {25}, Providence, R. I., 1969.


\item{[\Wir]}
{E. Wirsing},
{\it Approximation mit algebraischen Zahlen beschr\"ankten
Grades}, J. reine angew. Math. {206} (1961), 67--77.

\vskip1cm

\noindent Yann Bugeaud  \hfill{Michel Laurent}

\noindent Universit\'e Louis Pasteur
\hfill{Institut de Math\'ematiques de Luminy}

\noindent U. F. R. de math\'ematiques
\hfill{C.N.R.S. -  U.P.R. 9016 - case 907}

\noindent 7, rue Ren\'e Descartes      \hfill{163, avenue de Luminy}

\noindent 67084 STRASBOURG  (FRANCE)
\hfill{13288 MARSEILLE CEDEX 9  (FRANCE)}

\noindent {\tt bugeaud@math.u-strasbg.fr}
\hfill{{\tt laurent@iml.univ-mrs.fr}}

\end